\documentclass[11pt]{article}

\usepackage{amsmath,amsthm,amsfonts,amssymb,latexsym}
\usepackage{amsmath}

\usepackage{a4wide}

\newtheorem{thm}{Theorem}[section]

\newtheorem{lem}[thm]{Lemma}

\newtheorem{ex}[thm]{Example}

\begin{document}


\date{}\title{Deformation of $\mathfrak{sl}(2)$ and $\mathfrak{osp}(1|2)$-Modules of Symbols}

\author{Imed Basdouri\thanks{D\'epartement de Math\'ematiques, Facult\'e des Sciences de
Gafsa, Zarroug 2112 Gafsa, Tunisie. ~~~~~~~~~~~~~~ E.mail :
basdourimed@yahoo.fr}\and Mabrouk Ben Ammar\thanks{D\'epartement de
Math\'ematiques, Facult\'e des Sciences de Sfax, BP 802, 3038 Sfax,
Tunisie. ~~~~~~~~~~~~~~~~~~~~~~ E.mail :
mabrouk.benammar@fss.rnu.tn}}
\maketitle
\begin{abstract}
We consider the $\mathfrak{sl}(2)$-module structure on the spaces of
symbols of differential operators acting on the spaces of weighted
densities. We compute the necessary and sufficient integrability
conditions of a given infinitesimal deformation of this structure
and we prove that any formal deformation is equivalent to its
infinitesimal part. We study also the super analogue of this problem
getting the same results.
\end{abstract}

\maketitle {\bf Mathematics Subject Classification} (2010). 17B56,
53D55, 58H15.


\section{Introduction}
Let $\frak{vect}(1)$ be the Lie algebra of polynomial vector fields
on $\mathbb{R}$. Denote by $\mathcal{F}_\lambda=
\left\{fdx^{\lambda}\mid f\in \mathbb{R}[x]\right\}$ the space of
polynomial weighted densities of weight $\lambda\in\mathbb{R}$. The
space $\mathcal{F}_\lambda$ is a $\mathfrak{vect}(1)$-module for the
action defined by $$ L_{g\frac{d}{dx}}^\lambda(fdx^\lambda)=
(gf+\lambda g'f)dx^\lambda.$$

Any differential operator $A$ on $\mathbb{R}$ can be viewed  as the
linear mapping $f(dx)^\lambda\mapsto(Af)(dx)^\mu$ from
$\mathcal{F}_\lambda$ to $\mathcal{F}_\mu$ ($\lambda$, $\mu$ in
$\mathbb R$). Thus the space of differential operators is a
$\mathfrak{vect}(1)$-module, denoted $\mathrm
D_{\lambda,\mu}:=\mathrm{Hom}_\mathrm{diff}(\mathcal{F}_\lambda,
\mathcal{F}_\mu)$. The $\mathfrak{vect}(1)$ action is:
\begin{equation}\label{Lieder2}
L_X^{\lambda,\mu}(A)=L_X^\mu\circ A-A\circ L_X^\lambda.
\end{equation}
Each module $\mathrm{D}_{\lambda,\mu}$ has a natural filtration by
the order of differential operators; the graded module $\mathcal
S_{\lambda,\mu}:=gr\mathrm{D}_{\lambda,\mu}$ is called the
\textsl{space of symbols}. The quotient-module
$\mathrm{D}^k_{\lambda,\mu}/\mathrm{D}^{k-1}_{\lambda,\mu}$ is
isomorphic to module of tensor densities $\mathcal
F_{\lambda-\mu-k}$, the isomorphism is provided by the principal
symbol $\sigma_{pr}$ defined by
$$
A=\sum_{i=0}^ka_i(x)\partial^i_x\mapsto\sigma_{pr}(A)=a_i(x)(dx)^{\mu-\lambda-k}$$
As a $\frak{vect}(1)$-module, the space $\mathcal S_{\lambda,\mu}$
depends only on the difference $\delta=\mu-\lambda,$ so that
$\mathcal S_{\lambda,\mu}$ can be written as $\mathcal S_{\delta}$,
and we have
$$
\mathcal{S}_{\delta}=\bigoplus_{k=0}^{\infty}\mathcal{F}_{\delta-k}.$$
Denote by $\mathrm{D}_\delta$ the $\mathfrak{vect}(1)$-module of
differential operators in $\mathcal{S}_\delta$.

 The space $\mathrm{D}_{\lambda,\mu}$ cannot be
isomorphic as a $\mathfrak{vect}(1)$-module to the corresponding
space of symbols, but is a deformation of this space in the sense of
Richardson-Neijenhuis \cite{nr}.

If we restrict ourselves to the Lie subalgebra of
$\mathfrak{vect}(1)$ generated by
$\left\{\frac{d}{dx},\,x\frac{d}{dx},\,x^2\frac{d}{dx}\right\}$,
isomorphic to $\mathfrak{sl}(2)$, we get a family of infinite
dimensional $\mathfrak{sl}(2)$ modules, still denoted
$\mathcal{F}_{\lambda}$, $\mathrm{D}_{\lambda,\mu}$ and
$\mathrm{D}_\delta$.

We are also interested in the study of the analogue super
structures, namely, we consider the superspace $\mathbb{R}^{1|1}$
with coordinates $(x,\theta)$ where $\theta$ is the odd variable:
$\theta^2=0$. This superspace is equipped with the standard contact
structure given by the distribution $\langle\overline{\eta}\rangle$
generated by the vector field
$\overline{\eta}=\partial_{\theta}-\theta\partial_x$. That is, the
distribution $\langle\overline{\eta}\rangle$ is the kernel of the
following $1$-form:
\begin{equation*}
\alpha=dx+\theta d\theta.
\end{equation*}
Consider the superspace of polynomial functions
$$
\mathbb{R}[x,\theta]=\left\{F(x,\theta)=f_0(x)+\theta
f_1(x)~|~f_0,\,f_1\in\mathbb{R}[x]\right\}
$$ and consider the superspace
$\mathcal{K}(1)$ of contact polynomial vector fields on
$\mathbb{R}^{1|1}$. That is, $\mathcal{K}(1)$ is the superspace of
polynomial vector fields on $\mathbb{R}^{1|1}$ preserving the
distribution $\langle\overline{\eta}\rangle$:
$$
\mathcal{K}(1)=\big\{X\in\mathrm{Vect_{Pol}}(\mathbb{R}^{1|1})~|~[X,\,\overline{\eta}]=
F_X\overline{\eta}\quad\hbox{for some}~F_X\in
{\mathbb{R}}[x,\,\theta]\big\}. $$ The Lie superalgebra
$\mathcal{K}(1)$ is spanned by the vector fields of the form:
\begin{equation*}
X_F=F\partial_x
-\frac{1}{2}(-1)^{p(F)}\overline{\eta}(F)\overline{\eta},\quad\text{where}\quad
F\in \mathbb{R}[x,\theta].
\end{equation*}
We introduce the superspace
$\mathfrak{F}_\lambda=\{F\alpha^\lambda~|~F\in\mathbb{R}[x,\theta]\}$
of $\lambda$-densities on $\mathbb{R}^{1|1}$. This space is a
$\mathcal{K}(1)$-module for the action defined by
$$
\mathfrak{L}_{X_G}^\lambda(F\alpha^\lambda)=(X_G+\lambda G')(F)\alpha^\lambda.
$$
Similarly, we consider the $\mathcal{K}(1)$-module of linear
differential operators, $\mathfrak{D}_{\nu,\mu} : =
\mathrm{Hom}_\mathrm{diff} (\mathfrak{F}_{\nu},\mathfrak{F}_{\mu})$,
which is the super analogue of the space $\mathrm{D}_{\nu,\mu}$. The
$\mathcal{K}(1)$-action on $\mathfrak{D}_{\nu,\mu}$ is given by
\begin{equation}\label{d-action}
\mathfrak{L}^{\lambda,\mu}_{X_F}(A)=\mathfrak{L}^{\mu}_{X_F}\circ
A-(-1)^{p(A)p(F)} A\circ \mathfrak{L}^{\lambda}_{X_F}.
\end{equation}
The Lie superalgebra $\mathfrak{osp}(1|2)$, a super analogue of
$\mathrm{\frak {sl}}(2)$, can be realized as a subalgebra of
$\mathcal{K}(1)$:\begin{equation*}\mathfrak{osp}(1|2)=\text{Span}\left(X_1,\,X_{\theta},\,
X_{x},\,X_{x\theta},\, X_{x^2}\right).\end{equation*} The space of
even elements of $\mathfrak{osp}(1|2)$ is isomorphic to
$\mathfrak{sl}(2)$:
\begin{equation*}
(\mathfrak{osp}(1|2))_0=\text{Span}(X_1,\,X_{x},\,X_{x^2})=\mathfrak{sl}(2).
\end{equation*}
The super analogue of the space ${\cal S}_\delta$ is naturally the
superspace (see \cite{gmo}):
\begin{equation*}
\mathfrak{S}_{\delta}=\bigoplus_{k\in
\mathbb{N}}\mathfrak{F}_{\delta-\frac{k}{2}}.
\end{equation*}
Denote by $\mathfrak{D}_\delta$  the $\mathcal{K}(1)$-module of
linear differential operators in $\mathfrak{S_\delta}$.

In this paper, we study the deformations of the structure of the
$\mathfrak{sl}(2)$-modules $\mathcal{S}_\delta$ and their analogues
the $\mathfrak{osp}(1|2)$-modules $\mathfrak{S}_\delta$. We exhibit
the necessary and sufficient integrability conditions of a given
infinitesimal deformation. We prove that any formal deformation is
equivalent to its infinitesimal part and we give an example of
deformation with one parameter.

\section{Deformation }
Deformation theory of Lie algebra homomorphisms was first considered
with only one-parameter of deformation \cite{ff, nr, r}. Recently,
deformations of Lie (super)algebras with multi-parameters were
intensively studied ( see, e.g., \cite{aalo, abbo, bbbbk, bbdo, bb,
or, or1, or2}).

Let $\rho_0 :\mathfrak g \rightarrow\mathrm{End}(V)$ be an action of
a Lie (super)algebra $\mathfrak g$ on a vector (super)space $V$. It
is well known that the first cohomology space $\mathrm{H}^1(\frak
g;\mathrm{End}(V))$ determines and classifies infinitesimal
deformations up to equivalence. Thus, if
$\dim{\mathrm{H}^1(\mathfrak g;\mathrm{End}(V))}=m$, then choose
1-cocycles $\Upsilon_1,\ldots,\Upsilon_m$ representing a basis of
$\mathrm{H}^1(\mathfrak g;{\rm End}(V))$ and consider the
infinitesimal deformation
\begin{equation*}
 \rho=\rho_0+\sum_{i=1}^m{}t_i\,\Upsilon_i,
\end{equation*}
where $t_1,\ldots,t_m$ are independent parameters with
$p(t_i)=p(\Upsilon_i)$. We try to extend this infinitesimal
deformation to a formal one:
\begin{equation*}
 \rho= \rho_0+\sum_{i=1}^m{}t_i\,\Upsilon_i+
\sum_{i,j}{}t_it_j\, \rho^{(2)}_{ij}+\cdots,
\end{equation*}
where $\rho^{(2)}_{ij},\rho^{(3)}_{ijk},\ldots$ are linear maps from
$\frak g$ to $\mathrm{End}(V)$ with $p(\rho^{(2)}_{ij})=p(t_it_j),
p(\rho^{(3)}_{ijk})=p(t_it_jt_k),\dots$ such that
\begin{equation}\label{crochet}[\rho(x),\rho(y)=\rho([x,y]),\quad
 x,\,y\in\mathfrak{g}.\end{equation}
All the obstructions become from the condition (\ref{crochet}) and
it is well known that they lie in  $\mathrm{H}^2 (\frak g ,
\mathrm{End}(V))$. Thus, we will impose extra algebraic relations on
the parameters $t_1,\ldots,t_m$. Let $\mathcal{R}$ be an ideal in
$\mathbb{C}[[t_1,\ldots,t_m]]$ generated by some set of relations,
the quotient
\begin{equation*}
\mathcal{A}=\mathbb{C}[[t_1,\ldots,t_m]]/\mathcal{R}
\end{equation*}
is a (super)commutative associative (super)algebra with unity, and
we can speak about  deformations with base $\mathcal{A}$, (see
\cite{bbbbk,ff} for details).

\subsection{Deformation of the $\mathfrak{sl}(2)$-Modules of Symbols}
Now we study the formal deformations of the $\frak{ sl}(2)$-module
structure on the space of symbols:
\begin{equation*}
\mathcal{S}_\delta=\bigoplus_{k\geq0}\mathcal{F}_{\delta-k}.
\end{equation*}
The infinitesimal deformations are described by the cohomology space
$$
\mathrm{H}^1_\mathrm{diff}(\mathfrak{sl}(2),\mathrm{D}_\delta)=
\bigoplus_{i,
j\geq0}\mathrm{H}^1_\mathrm{diff}\left(\mathfrak{sl}(2),
\mathrm{D}_{\delta-j,\delta-i}\right)
$$
where $\mathrm{H}^*_\mathrm{diff}$ denotes the differential
cohomology; that is, only cochains given by differential operators
are considered. In fact, Lecomte computed  $\mathrm{H}^1_{\rm
diff}\left(\mathfrak{sl}(2), \mathrm{D}_{\lambda,\lambda'}\right)$,
see \cite{lec}. He showed that non-zero cohomology
$\mathrm{H}^1_\mathrm{diff}\left(\mathfrak{vect}(1),
\mathrm{D}_{\lambda,\lambda'}\right)$ only appear if
$\lambda=\lambda'$ or
$(\lambda,\lambda')=(\frac{1-k}{2},\frac{1+k}{2})$ where
$k\in\mathbb{N}^*$. Thus, we distinguish two cases:
\begin{itemize}
  \item [(i)] If $\delta\notin\frac{1}{2}(\mathbb{N}+2)$, then
  $$\mathrm{H}^1_\mathrm{diff}(\mathfrak{sl}(2),\mathrm{D}_\delta)=\bigoplus_{k\geq0}
  \mathrm{H}^1_\mathrm{diff}\left(\mathfrak{sl}(2),
\mathrm{D}_{\delta-k,\delta-k}\right).$$ The space
$\mathrm{H}^1_\mathrm{diff}\left(\mathfrak{sl}(2),
\mathrm{D}_{\lambda,\lambda}\right)$ is one dimensional and it is
spanned by the cohomology classe of the cocycle $A_\lambda$ given by
\begin{equation*}\label{A}\begin{array}{l}A_\lambda(F\frac{d}{dx})(f{dx}^{\lambda})=
F'f{dx}^{\lambda}.\end{array}\end{equation*}
  \item [(ii)] If $2\delta=m\in(\mathbb{N}+2)$, then
\begin{equation*}
\mathrm{H}^1_\mathrm{diff}(\mathfrak{sl}(2),\mathrm{D}_\delta)=
\bigoplus_{k=[\frac{m+1}{2}]}^{m-1}\mathrm{H}^1_\mathrm{diff}
\left(\mathfrak{sl}(2),\mathrm{D}_{{\frac{m-2k}{2}},{\frac{2+2k-m}{2}}}
\right)\oplus\bigoplus_{k\geq0}\mathrm{H}^1_\mathrm{diff}\left(\mathfrak{sl}(2),
\mathrm{D}_{{\frac{m}{2}-k},{\frac{m}{2}}-k}\right).
\end{equation*}
The space ${\mathrm H}^1_\mathrm{diff}\left(\mathfrak{sl}(2),
\mathrm{D}_{{\frac{m-2k}{2}},{\frac{2+2k-m}{2}}}\right)$ is two
dimensional and spanned by the cohomology classes of the 1-cocycles,
$B_k$ and ${C}_k$ given by
\begin{align*}\begin{array}{llll}
B_k(F\frac{d}{dx})(f{dx}^{\frac{m-2k}{2}})
=F'f^{(2k-m+1)}{dx}^{\frac{2+2k-m}{2}},\,{C}_k(F\frac{d}{dx})(f{dx}^{\frac{m-2k}{2}})=
F''f^{(2k-m)}dx^{\frac{2+2k-m}{2}}.
\end{array}\end{align*}
\end{itemize}

In our study, an infinitesimal deformation of the
$\mathfrak{sl}(2)$-module structure on the space
$\mathcal{S}_\delta$ is of the form
\begin{equation}
\label{InDef}\mathcal{L}_X=L_X+\mathcal{L}^{(1)}_X,
\end{equation} where $L_X$ is the Lie derivative of $\mathrm{D}_\delta$
 along the vector field $X$ defined by
(\ref{Lieder2}), and
\begin{equation*}
\label{InfinDef2}
\mathcal{L}_X^{(1)}=\left\{\begin{array}{ll}\sum_{k\geq0}
a_k\,A_{\frac{m}{2}-k}(X)&\text{ if
}\delta\notin\frac{1}{2}(\mathbb{N}+2)\\[6pt]
\sum_{k\geq0}{}a_k\,A_{\frac{m}{2}-k}(X)+\sum_{k=[\frac{m+1}{2}]}^{m-1}
\big(b_k\,B_{k}(X)+c_{k}C_{k}(X)\big)&\text{ if
}2\delta=m\in(\mathbb{N}+2),\end{array}\right.
\end{equation*}
and where $a_{k}$, $b_{k}$  and $c_{k}$ are independent parameters.
\begin{thm}
\label{Main12}  The following  conditions are necessary and
sufficient  for integrability of the infinitesimal deformation
(\ref{InDef}):
\begin{equation}\label{condition}\begin{array}{l}(2k-m+1)b_{k}a_{m-k-1}+c_{k}a_{k}-c_{k}a_{m-k-1}=0,
\quad[\frac{m+1}{2}]\leq k\leq m-1\end{array}.\end{equation}
Moreover, any formal deformation is equivalent to its infinitesimal
part.
\end{thm}
\begin{proofname}. Note that if
$\delta\notin\frac{1}{2}(\mathbb{N}+2)$ then the parameters $b_k$
and $c_k$ can be assumed to be zero, and then, there are no
integrability conditions. Assume that the infinitesimal deformation
(\ref{InDef}) can be integrated to a formal deformation

\begin{equation*}
\label{formal1} \mathcal{L}_X=L_X+\mathcal{L}^{(1)}_X+
\mathcal{L}^{(2)}_X+\mathcal{L}^{(3)}_X+\cdots
\end{equation*}
where $\mathcal{L}^{(1)}_X$ is given by (\ref{InfinDef2}) and
$\mathcal{L}^{(2)}_X$ is a quadratic polynomial in $a_k$, $b_k$ and
$c_k$ with coefficients in $\mathcal{D}_{\delta}$. We compute the
conditions for the second-order terms $\mathcal{L}^{(2) }$. Consider
the quadratic terms of the homomorphism condition
\begin{equation}
\label{HomomCond2}
[\mathcal{L}_X,\mathcal{L}_Y]=\mathcal{L}_{[X,Y]}.
\end{equation}
By a straightforward computation, the homomorphism condition
(\ref{HomomCond2}) gives for the second-order terms the following
 equation
\begin{equation}\label{L2}
 \delta({\cal
L}^{(2)})=\frac{1}{2}\sum_{k=[\frac{m+1}{2}]}^{m-1}((2k-m+1)b_{k}a_{m-k-1}+c_{k}a_{k}-c_{k}a_{m-k-1})\Phi_{2k-m+1},
\end{equation}
where $\Phi_{k}$ is the nontrivial 2 cocycle given by
\begin{equation}\label{Phi}\begin{array}{llll}\Phi_k(F\frac{d}{dx},
G\frac{d}{dx})( f dx^\frac{1-k}{2})
=(F^{'}G^{''}-F^{''}G^{'})f^{(k-1)}dx^\frac{1+k}{2}\end{array}\end{equation}
The condition (\ref{condition}) is necessary since the operator
$\Phi_{k}$ is a nontrivial 2 cocycle  spanning the space
$\mathrm{H}^2_\mathrm{diff}\left(\mathfrak{sl}(2),
\mathrm{D}_{\frac{m-2k}{2},\frac{2+2k-m}{2}}\right)$, (see
\cite{lec}). The solution $\mathcal{L}^{(2)}$ of (\ref{L2}) can be
chosen identically zero. Choosing the highest-order terms
$\mathcal{L}^{(m)}$ with $m\geq3$, also identically zero, one
obviously obtains a deformation (which is of order 1 in $a_k$, $b_k$
and $c_k$). \hfill$\Box$
\end{proofname}

{\ex \label{Example1} Let us consider
$\delta=\frac{m}{2}\in\frac{1}{2}(\mathbb{N}+2)$ and let
$(\alpha_k)_{k\geq0}$ be sequence of real numbers such that, for
$[\frac{m+1}{2}]\leq k\leq m-1$, we have
$\alpha_k\neq\alpha_{m-k-1}$. Put $b_k=t$, $a_k=\alpha_k\, t$ and
$c_k=\frac{(2k-m+1)\alpha_k }{\alpha_k-\alpha_{m-k-1}}\,t$. So, we obtain a
deformation of $\mathcal{S}_\delta$ with one parameter $t$:
\begin{equation*}
\mathcal{L}=
L+t\sum_{k\geq0}{}\alpha_k\,A_{\frac{m}{2}-k}(X)+t\sum_{k=[\frac{m+1}{2}]}^{m-1}
\left(B_{k}(X)+\frac{(2k-m+1)\alpha_k}{\alpha_k-\alpha_{m-k-1}}\,C_{k}(X)\right).
\end{equation*}}
Of course it is easy to give many other examples of true
deformations with one parameter or with several parameters.
\subsection{Deformation of the $\mathfrak{osp}(1|2)$-Modules of Symbols}
We study the super analogous of the previous case. That is, we study
deformations of the $\mathfrak{osp}(1|2)$-module of differential
linear operators in the space of symbols on $\mathbb{R}^{1|1}$:
\begin{equation*}
\mathfrak{S}_\delta=\bigoplus_{k\geq0}\frak F_{\delta-\frac{k}{2}}.
\end{equation*}
The infinitesimal deformations are described by the cohomology space
$$
\mathrm{H}^1_\mathrm{diff}(\mathfrak{osp}(1|2),\frak{D}_\delta)=
\bigoplus_{i, j\geq0}\mathrm{H}^1_{\rm
diff}\left(\mathfrak{osp}(1|2),
\mathfrak{D}_{\delta-\frac{j}{2},\delta-\frac{i}{2}}\right).
$$
In \cite{bb1}, it was proved that non-zero cohomology
$\mathrm{H}^1_{\rm
diff}\left(\mathcal{K}(1),\frak{D}_{\lambda,\lambda'}\right)$ only
appear if $\lambda=\lambda'$ or
$(\lambda,\lambda')=(\frac{1-k}{2},\frac{k}{2})$ where $k\in
\mathbb{N}$. Thus, as before, we have to distinguish two cases:
\begin{itemize}
  \item [(i)] If $\delta\notin\frac{1}{2}(\mathbb{N}+1)$, then
  $$\mathrm{H}^1_\mathrm{diff}(\mathfrak{osp}(1|2),\frak{D}_\delta)=\bigoplus_{k\geq0}
  \mathrm{H}^1_{\rm diff}\left(\mathfrak{osp}(1|2),
\frak{D}_{\delta-\frac{k}{2},\delta-\frac{k}{2}}\right).$$ The space
$\mathrm{H}^1_\mathrm{diff}\left(\mathfrak{osp}(1|2),
\frak{D}_{\frac{2\delta-k}{2},\frac{2\delta-k}{2}}\right)$ is one
dimensional and it is spanned by the cohomology classe of the
cocycle $\Upsilon'_{2\delta-k}$ given by
\begin{equation*}\begin{array}{l}\Upsilon'_{2\delta-k}(F\frac{d}{dx})=
F'.\end{array}\end{equation*}
  \item [(ii)] If $2\delta=m\in(\mathbb{N}+1)$, then
\begin{equation*}
\mathrm{H}^1_\mathrm{diff}(\mathfrak{osp}(1|2),\mathfrak{D}_\delta)=
\bigoplus_{k=1}^{m}\mathrm{H}^1_{\rm diff}\left(\mathfrak{osp}(1|2),
\mathfrak{D}_{{\frac{1-k}{2}},{\frac{k}{2}}}\right)\oplus\bigoplus_{k\geq0}\mathrm{H}^1_\mathrm{
diff}\left(\mathfrak{osp}(1|2),
\mathfrak{D}_{{\frac{m-k}{2}},{\frac{m-k}{2}}}\right).
\end{equation*}
The space $\mathrm{H}^1_\mathrm{diff}\left({\rm\frak
{osp}}(1|2),\mathfrak{D}_{{\frac{1-k}{2}},{\frac{k}{2}}}\right)$ is
two dimensional and spanned by the cohomology classes of the
1-cocycles, $\Upsilon_k$ and $\widetilde{\Upsilon}_k$ given by
\begin{equation*}\label{uk}\begin{array}{ll}
\Upsilon_k(X_G)=(-1)^{|G|}{\eta}^2(G)\overline{\eta}^{2k-1}, \quad
\widetilde{\Upsilon}_k(X_G)=(-1)^{|G|}((k-1)\eta^4(G)\overline{\eta}^{2k-3}
+\eta^3(G)\overline{\eta}^{2k-2}).\end{array}
\end{equation*}
\end{itemize}

Any infinitesimal  deformation of the $\mathfrak{spo}(1|2)$-module
structure on  $ \mathfrak{S}_\delta$  is of the form
\begin{equation}\label{infdef2}
\widetilde{\mathfrak
{L}}_{X_F}=\mathfrak{L}_{X_F}+\mathfrak{L}^{(1)}_{X_F},
\end{equation}where $\mathfrak{L}_{X_F}$ is the Lie derivative of $\mathfrak{D}_\delta$
 along the vector field $X_F$ defined by
(\ref{d-action}), and
\begin{equation*}
 \mathfrak{L}_{X_F}^{(1)}=\left\{\begin{array}{ll}\sum_{k\geq0}
\mathfrak{a}_{2\delta-k}\,\Upsilon'_{2\delta-k}(X_F)&\text{ if
}\delta\notin(\frac{1}{2}\mathbb{N}+1)\\[6pt]
\sum_{k\geq0}{}\mathfrak{a}_{m-k}\,\Upsilon'_{m-k}(X_F)+\sum_{k=1}^{m}
\big(\mathfrak{b}_k\,\Upsilon_{k}(X_F)+\mathfrak{c}_{k}\widetilde{\Upsilon}_{k}(X_F)\big)&\text{
if }2\delta=m\in(\mathbb{N}+1),\end{array}\right.
\end{equation*}
and where $\mathfrak{a}_{k}$, $\mathfrak{b}_{k}$  and
$\mathfrak{c}_{k}$ are independent parameters.

Our main result in the super setting is the following
\begin{thm} The following  conditions are necessary and
sufficient  for integrability of the infinitesimal deformation
(\ref{infdef2}):
\label{Main22}\begin{equation}\label{condition2}\begin{array}{l}\frak{b}_{k}\frak{a}_{1-k}-
\frak{c}_{k}\frak{a}_{k}+\frak{c}_{k}\frak{a}_{1-k}=0, \quad 1\leq
k\leq m\end{array}.\end{equation} Moreover, any formal deformation
is equivalent to its infinitesimal part.
\end{thm}

\begin{proofname}. Assume that the infinitesimal deformation
(\ref{infdef2}) can be integrated to a formal deformation:
\begin{equation*}
\widetilde{\mathfrak L}_{X_F}=\frak{L}_{X_F}+{\frak
L}^{(1)}_{X_F}+\mathfrak{L}^{(2)}_{X_F}+\cdots
\end{equation*}
By a straightforward computation, the homomorphism condition
\begin{equation*}
 [\widetilde{\frak L}_{X_F},\widetilde{\mathfrak
L}_{X_G}]=\widetilde{\mathfrak L}_{X_{\{F,G\}}}
\end{equation*} gives for the second-order terms the following
 equation
$$\delta(\mathfrak{L}^{(2)})=\sum_{k=1}^m(\mathfrak{b}_{k}\mathfrak{a}_{1-k}-
\mathfrak{c}_{k}\mathfrak{a}_{k}+\mathfrak{c}_{k}\mathfrak{a}_{1-k})\Omega_k$$
where
$\Omega_{k}:\mathfrak{osp}(1|2)\times\mathfrak{osp}(1|2)\rightarrow\mathfrak{D}_{\frac{1-k}{2},\frac{k}{2}}$
is defined by
$$
\Omega_k(X_F,X_G)=(-1)^{p(F)+p(G)}(k-1)(F'G''-F''G')\overline{\eta}^{2k-3}+(\overline{\eta}(F')G'-
(-1)^{p(F)p(G)}F'\overline{\eta}(G'))\overline{\eta}^{2k-2}.
$$
We ill prove the following lemma and then we conclude as for Theorem
\ref{Main12}.
\begin{lem}\label{cup3} The map $\Omega_k$ is a nontrivial odd
2 cocycle.
\end{lem}
\begin{proofname}.
The map $\Omega_k$ is a 2 cocycle since it is the cup-product of 1
cocycles. It is easy to see that $\Omega_k$ is an odd map, so,
$\Omega_k(\mathfrak{sl}(2)\times\mathfrak{sl}(2))\subset
\left(\mathfrak{D}_{\frac{1-k}{2},\frac{k}{2}}\right)_1$. In
\cite{bbbbk} it was proved that, as $\mathfrak{sl}(2)$-module, we
have
\begin{equation}\label{decomp}\left(\mathfrak{D}_{\frac{1-k}{2},\frac{k}{2}}\right)_1\simeq
\Pi\left(\mathrm{D}_{\frac{2-k}{2},\frac{k}{2}}\oplus\mathrm{D}_{\frac{1-k}{2},\frac{1+k}{2}}\right)\end{equation}
where $\Pi$ is the change of parity operator. We check that the
restriction of $\Omega_k$ to
$\mathfrak{sl}(2)\times\mathfrak{sl}(2)$ is a nontrivial 2 cocycle.
Indeed, let $X_F,\,
X_G\in\mathfrak{sl}(2)\subset\mathfrak{osp}(1|2)$, it is easy to see
that
$$
(-1)^k\Omega_k(X_F,X_G)=(k-1)\Phi_{k-1}(X_F,X_G)
\circ\partial_\theta-k\theta\Phi_{k}(X_F,X_G),
$$
or equivalently, according to the decomposition (\ref{decomp}), we
have
$$
(-1)^k\Omega_k|_{\mathfrak{sl}(2)\times\mathfrak{sl}(2)}=\Pi\circ((k-1)\Phi_{k-1}-k\Phi_{k})
$$
where $\Phi_k$ is the nontrivial 2 cocycle defined by (\ref{Phi}).
Thus, $\Omega$ is a nontrivial 2 cocycle. \hfill$\Box$
\end{proofname}
\end{proofname}
 Obviously, as for the $\mathfrak{sl}(2)$-module
 $\mathcal{S}_\delta$, it easy to construct many examples of true
deformations of the $\mathfrak{osp}(1|2)$-module
$\mathfrak{S}_\delta$ with one parameter or with several parameters.


\begin{thebibliography}{99}

\bibitem{aalo}
B. Agrebaoui, F. Ammar, P. Lecomte, V. Ovsienko, {\it
Multi-parameter deformations of the module of symbols of
differential operators}, Internat. Mathem. Research Notices, 2002,
N¡. 16, 847--869.
\bibitem{abbo}B. Agrebaoui, N. Ben Fraj, M. Ben Ammar, and V.
Ovsienko, {\it Deformation of modules of differential forms,}
NonLinear Mathematical Physics, vol. 10(2003)num. 2, 148--156.


\bibitem{bb1}I. Basdouri, M. Ben Ammar, {\it Cohomology of $\frak {osp}(1|2)$
Acting on Linear Differential Operators on the Supercercle $S^{1|1}$.}
 Letters in Mathematical Physics(2007)81:239--251.


\bibitem{bbbbk}
I. Basdouri, M. Ben Ammar, N. Ben Fraj, M. Boujelbene and K.
Kammoun, {\it Cohomology of the Lie Superalgebra of Contact Vector
Fields on $\mathbb{R}^{1|1} $ and Deformations of the Superspace of
Symbols .} Jour of Nonlinear Math Physics, Vol. 16, No. 4 (2009)
1--37.


\bibitem{bbdo}I. Basdouri, M. Ben Ammar, B. Dali and S. Omri,
{\it Deformation of $\mathfrak{vect}(1)$-Modules of Symbols.}
Journal of Geometry and Physics 60 (2010) 531--542


\bibitem{bb} M. Ben Ammar, M. Boujelbene,
{\it $\mathfrak{sl}(2)$-Trivial Deformation of
$\mathrm{Vect_{Pol}}(\mathbb{R})$-Modules of Symbols.} SIGMA 4
(2008), 065, 19 pages.

\bibitem{ff}
A. Fialowski, D. B. Fuchs, {\it Construction of miniversal
deformations of Lie algebras} J. Func. Anal. {\bf 161:1} (1999)
76--110.



\bibitem{gmo} H. Gargoubi, N. Mellouli and V. Ovsienko
{\it Differential Operators on Supercircle: Conformally
Equivariant Quantization and Symbol Calculus}, Letters
in Mathematical Physics (2007) {\bf 79}: 51–65.



\bibitem{lec}P. B. A. Lecomte, {\it On the cohomology of
$\frak{sl}(n + 1;\mathbb{R})$ acting on differential operators and
$\frak{sl}(n + 1;\mathbb{R})$-equivariant symbols,} Indag. Math. NS.
11 (1), (2000), 95 114.

\bibitem{nr}
A. Nijenuis, R. W. Richardson Jr., {\it Deformations of
homomorphisms of Lie groups and Lie algebras}, Bull. Amer. Math.
Soc. {\bf 73} (1967), 175--179.

\bibitem{or}
V. Ovsienko,  C. Roger, {\it Deforming the Lie algebra of vector
fields on~$S^1$ inside the Lie algebra of pseudodifferential
operators on $S^1$}, AMS Transl. Ser.~2, (Adv. Math. Sci.) vol.~194
(1999) 211--227.

\bibitem{or1}
 V. Ovsienko, C. Roger,
{\it Deforming the Lie algebra of vector fields on~$S^1$ inside the
Poisson algebra on $\dot T^*S^1$}, Comm. Math. Phys., {\bf 198}
(1998) 97--110.

\bibitem{or2}
 V. Ovsienko, C. Roger,
{\it Deforming the Lie algebra of vector fields on~$S^1$ inside the
Poisson algebra on $\dot T^*S^1$}, Comm. Math. Phys., {\bf 198}
(1998) 97--110.

\bibitem{r} R.W. Richardson, {\it Deformations of subalgebras of Lie algebras}, J. Diff. Geom. , 3,
(1969), 289–308.

\end{thebibliography}
\end{document}